\documentclass[12pt,a4paper]{article}
\usepackage[utf8]{inputenc}
\usepackage{amsmath,amsthm}
\usepackage{amsfonts,graphicx}
\usepackage{amssymb,hyperref}
\usepackage{color,url}
\usepackage{natbib}      
\usepackage{authblk}

\newtheorem{defn}{Definition}

\begin{document}

\title{An $r$-convex set which is not locally contractible}

\author[1]{Alejandro Cholaquidis}
 
\affil[1]{Facultad de Ciencias - Universidad de la República. Uruguay}

 \maketitle

\begin{abstract}
	The study of shape restrictions of subsets of $\mathbb{R}^d$ have several applications in many areas, being convexity, $r$-convexity, and positive reach, some  of the most famous, and typically imposed in set estimation. The following problem was attributed to K. Borsuk, by J. Perkal in 1956: find an $r$-convex set which is not locally contractible. Stated in that way is trivial to find such a set.   However, if we ask the set to be equal to the closure of its interior (a condition fulfilled for instance if the set is the support of a probability distribution absolutely continuous with respect to the $d$-dimensional Lebesgue measure), the problem is much more difficult. We present a counter example of a not-locally contractible set, which is $r$-convex. This also proves that the class of supports with positive reach of absolutely continuous distributions includes strictly the class of $r$-convex supports. 
\end{abstract}

%% Keywords should be separated by \*\ sign
\textit{keywords} $r$-convex set \*\ locally contractible set. \*\  positive reach

Several shape restriction have been introduced  to subsets of $\mathbb{R}^d$, being the convexity the most  studied. However, in many practical applications to assume that the set is convex is too restrictive. To overcome this limitation, and as a natural generalization, the notion of $r$-convex set was introduced. Let us recall that, given $r>0$, a set $S\subset \mathbb{R}^d$ is said to be $r$-convex, if $S=C_r(S),$ where $C_r(S)$ is the $r$-convex hull of $S$, i.e: the intersection of the complement of all open balls of radius $r$, not meeting $S$. More precisely,

\begin{equation} \label{rhull}
	C_r(S)=\bigcap_{\big\{ \mathring{B}(x,\alpha):\ \mathring{B}(x,r)\cap S=\emptyset\big\}} \Big(\mathring{B}(x,r)\Big)^c,
\end{equation}
where $\mathring{B}(x,r)$ is the open ball of radius $r$ centered at $x$. Convex sets are $r$-convex, for any $r>0$, the converse implication being trivially false. The notion of $r$-convex set has been studied in depth in set estimation, several results were obtained by assuming that the set is $r$-convex, see for instance \cite{pat08}, \cite{rod07}, \cite{wal97}, \cite{cue10}.

Another well known shape restriction, introduced by Federer (see \cite{fed59}) is the notion of set with positive reach. It is defined as the largest distance from which any point outside $S$ has a unique nearest point in $S$. In \cite{cue12} it is proven that the class of sets with reach $r$ contains the class of $r$-convex sets and the inclusion is strict. Federer proved that the sets with positive reach are locally contractible. In \cite{perk} posed the question (attributed to K. Borsuk) if that is also true for the class of $r$-convex sets. Stated in that way the assertion is trivially false as we will see. 
However,   imposing the condition that the set equals the closure of its interior (which is fulfilled if the set is the support of a probability distribution absolutely w.r.t the $d$-dimensional Lebesgue measure) the  problem becomes non-trivial.  We prove that the assertion is also false, with this additional restriction. The counter example proves then that the class of $r$-convex sets contains strictly the class of sets with positive reach. Let us recall the notions of locally contractible.

\begin{defn} \label{loc-contr} A topological space $X$ is said to be weakly locally contractible at a point $x$ if for every neighborhood $U_x$ of $x$ there exists $V_x\subset U_x$ a neighborhood such that $V_x$ is homotopy equivalent in $U_x$ (that is, with respect to the subspace topology) to a one-point space. A topological space $X$ is weakly locally contractible  if it is locally contractible at every point. If every point has a local base of contractible neighbourhoods its is said to be strongly locally contractible.
\end{defn}

To find subsets of $\mathbb{R}^d$ fulfilling the $r$-convexity property for some $r>0$ but not strongly locally contractible is quite easy if the condition $S=\overline{\mathring{S}}$ (i.e, the set is equal to the closure of its interior) is not required. We will present one later. Observe that this condition is not required in \cite{perk}, but it is quite important because, for instance, is fulfilled by the support of an absolutely continuous probability distribution (with respect to Lebesgue measure) in $\mathbb{R}^d$.  We think that this condition was wrongly omitted in \cite{perk} and also in \cite{ml}. We aim to present a compact $r$-convex subset of $\mathbb{R}^2$, not strongly contractible. \\

\section{The set}
Let $a_n=1/2^n$. Let us consider the following family of intervals of $\mathbb{R}$, $B_0=[0,a_1]$, $A_1=[a_1,a_1+a_2]$, $B_1=[a_1+a_2,a_1+2a_2]$, $A_2=[a_1+2a_2,a_1+2a_2+a_3]$, $B_2=[a_1+2a_2+a_3,a_1+2a_2+2a_3]$, in general 

$$A_n=[a_1+\sum_{j=2}^n 2a_j,a_1+\sum_{j=2}^n 2a_j+a_{n+1}]\quad \text{and} \quad B_n=[a_1+\sum_{j=2}^n 2a_j+a_{n+1},a_1+\sum_{j=2}^{n+1} 2a_j].$$
Observe that,
$$B_0\cup \bigcup_{n=1}^\infty A_n\cup \bigcup_{n=1}^\infty B_n\cup\{2\}\subset[0,2]$$
Let $\mathcal{C}=B_0\cup 	\bigcup_{n=1}^\infty B_n\cup\{2\}$. It is clear that $\{(x,0):x\in \mathcal{C}\}$ is compact, $r$-convex for any $r>0$ but it is not locally contractible, since $(2,0)$ is not locally connected. \\

\begin{figure}
	\centering
	\includegraphics[width=1\linewidth]{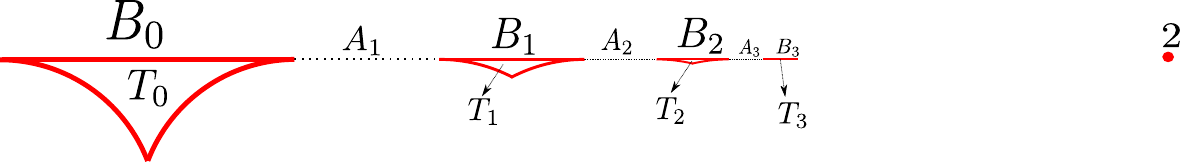}
	\caption{The set $\mathcal{T}=\cup_n T_n$ is compact, $r$-convex but not locally contractible at $(2,0)$}
	\label{fig:tn}
\end{figure}

Let us denote for $n\in \mathbb{N}$, $[b_n,b_{n+1}]$ the extremes of the intervals $B_n$. Let us consider $D_n= \mathring{B}[(b_n,-r),r]$ the open ball in $\mathbb{R}^2$ centred at $(b_n,-r)$. We denote by $T_n=\{z=(x,y): x\in B_n, y\leq 0, \text{ and } z\in (D_n\cup D_{n+1})^C\}$, see Figure \ref{fig:tn}.
The set
$$\mathcal{T}=(2,0)\cup\bigcup_{n=0}^\infty T_n$$
is compact, $r$-convex (this is immediate by construction), fulfills $\overline{\mathring{\mathcal{T}}}=\mathcal{T}$, but it is not weakly locally contractible at $(2,0)$ since any open set containing $(2,0)$ is not connected.

\end{document}